\documentclass[12pt]{article}

\usepackage{fullpage}
\usepackage{amsmath}
\usepackage{amssymb}
\usepackage{color}
\usepackage{graphicx}
\usepackage{authblk}

\usepackage[margin=0.6in]{geometry}

\title{How big a table do you need for your jigsaw puzzle?}

\author[1]{Madeleine Bonsma-Fisher}
\author[2]{Kent Bonsma-Fisher}
\date{}

\affil[1]{m.bonsma@utoronto.ca}
\affil[2]{kentfisher@gmail.com} 

\begin{document}
\maketitle

\section{Abstract}

Jigsaw puzzles are typically labeled with their finished area and number of pieces. With this information, is it possible to estimate the area required to lay each piece flat before assembly? We derive a simple formula based on two-dimensional circular packing and show that the unassembled puzzle area is $\sqrt{3}$ times the assembled puzzle area, independent of the number of pieces. We perform measurements on 9 puzzles ranging from $333 \text{ cm}^2$ (9 pieces) to $6798 \text{ cm}^2$ (2000 pieces) and show that the formula accurately predicts realistic assembly scenarios.

\section{Introduction}

Jigsaw puzzle enthusiasts know that the area required to lay all puzzle pieces flat is larger than the assembled puzzle area. To the best of our knowledge, no formula exists to accurately estimate the unassembled area of a puzzle. Such a formula would be helpful to assess whether the intended puzzle table can hold all the unassembled pieces in a single layer, the ideal puzzle setup.

Using principles from condensed matter physics \cite{Ashcroft76}, we developed a simple model approximating each puzzle piece as a square, then assuming that the unassembled pieces will pack together as if the squares were circumscribed by a circle. Using ideal 2D circular packing \cite{Chang2010}, we calculated the total area of the unassembled pieces and found that it depends on the area of the assembled puzzle with a very simple relationship: the unassembled area is $\sqrt{3}$ times the assembled area.

\section{Model}

\begin{figure}
\centering
\includegraphics[width=0.3\linewidth]{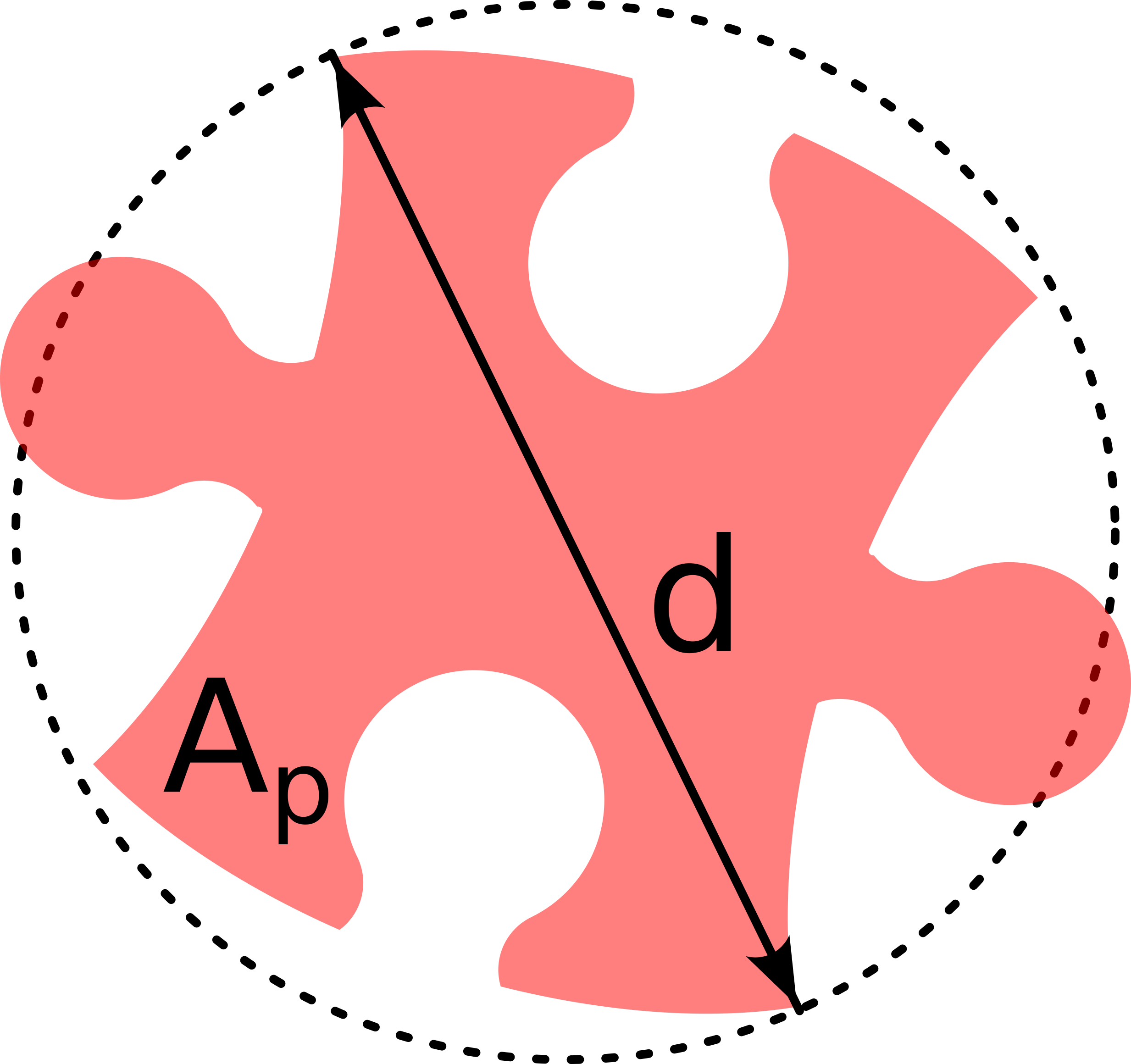} 
\caption{Diagram of a puzzle piece that is circumscribed by a circle of radius $d$. We assume $d$ to be equal to the diagonal of a square with area $A_p = A_a/N$, where $A_a$ is the assembled puzzle area and $N$ is the total number of pieces.} \label{fig:puzzle_diagram}
\end{figure}

Let $N$ be total number of puzzle pieces, $A_a$ be the area of the assembled puzzle, and $A_s$ be the area of the unassembled puzzle. We approximate each piece as a square of area $A_p = \frac{A_a}{N}$. Each puzzle piece edge has length $\sqrt{A_p} = \sqrt{A_a/N}$. Since before assembly the pieces are randomly oriented, we assume on average they behave like circles with a diameter equal to the square's diagonal (Figure \ref{fig:puzzle_diagram}).

\begin{figure}
\centering
\includegraphics[width=0.4\linewidth]{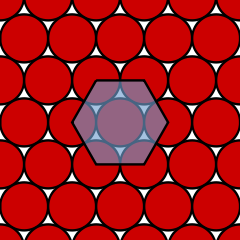} 
\caption{Hexagonal packing of circles on a 2D surface. Figure reproduced from Wikipedia, public domain.} \label{fig:hexagonal_lattice}
\end{figure}

We calculate the circle's diameter $d$ using the Pythagorean theorem, giving $d = \sqrt{2 A_p} = \sqrt{2 A_a/N}$. Next, we calculate the total unassembled area assuming the pieces pack on a hexagonal lattice (Figure \ref{fig:hexagonal_lattice}). The edge length of the hexagon is $d$, the circle diameter, and the area of a hexagon is $A_h = \frac{3\sqrt{3}}{2}d^2$. There are three full circles inside one hexagonal element: one in the centre, and each of the six edge circles contains $1/3$ of the area of a circle. The packing area of a puzzle piece is then 1/3 the hexagon area. 

For $N$ pieces with approximate area $A_h/3$ and hexagons with edge length equal to the diagonal of the approximate square puzzle piece $d = \sqrt{2 A_a/N}$, we arrive at equation \ref{eqn:area} for the area of the unassembled puzzle.

\begin{equation} \label{eqn:area}
A_s = N\frac{A_h}{3} = N\frac{\sqrt{3}}{2}d^2 = N\frac{\sqrt{3}}{2}\frac{2 A_a}{N} = \sqrt{3}A_a \approx 1.73 A_a
\end{equation}

Interestingly, the area of the unassembled puzzle is simply $\sqrt{3}$ times the area of the assembled puzzle, independent of the number of pieces.

\section{Results}

\begin{figure}
\centering
\includegraphics[width=0.7\linewidth]{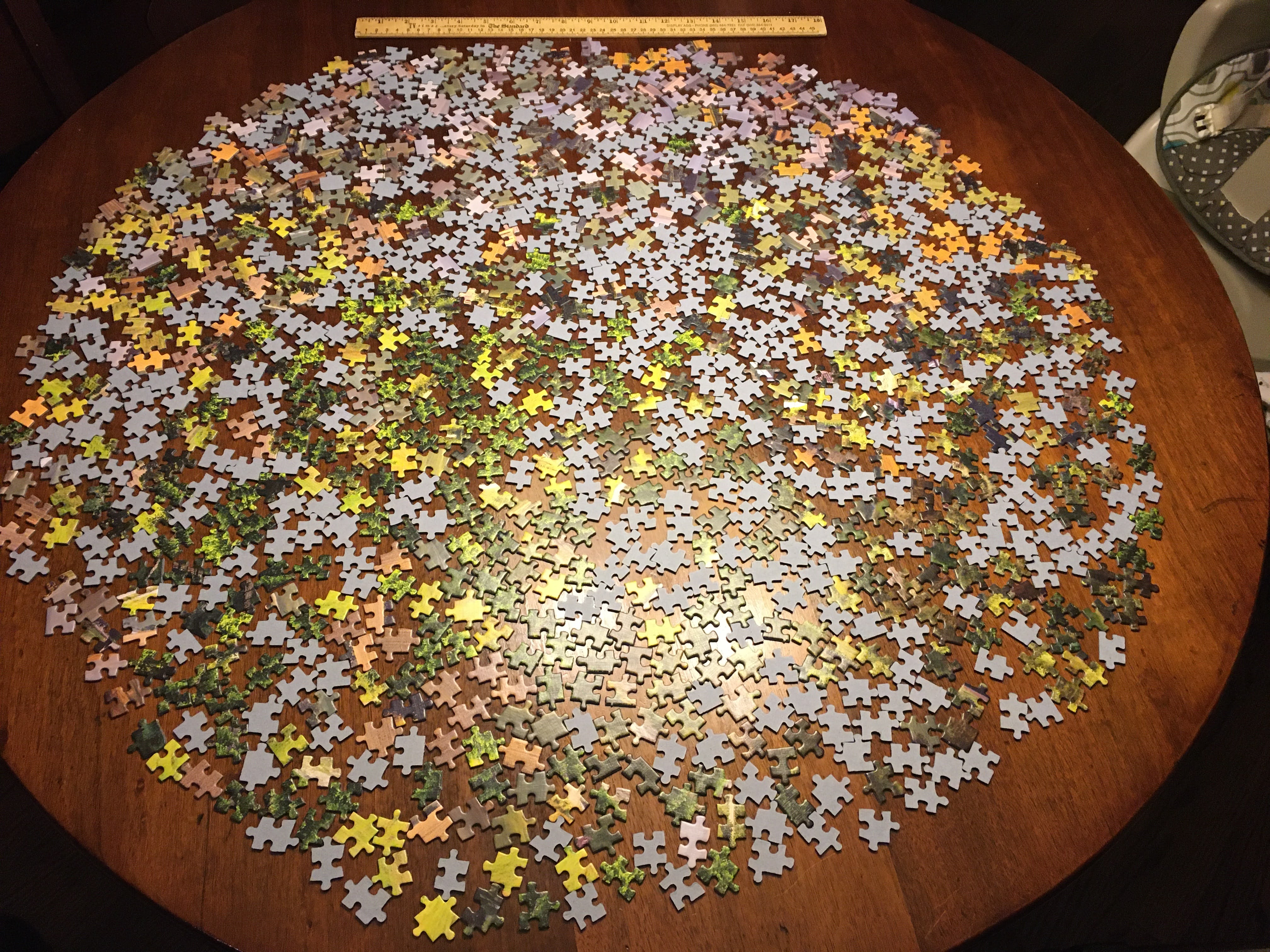} 
\caption{Photograph of an unassembled 1008-piece puzzle with pieces arranged in an approximate circular shape. We measured the long and short axes ($X_s$ and $Y_s$) of the resulting ellipse to estimate the unassembled area.} \label{fig:puzzle_measurement}
\end{figure}

We completed 9 different puzzles of varying sizes and numbers of pieces and measured their unassembled and assembled areas. Before assembly, we laid out all the pieces in a flat layer in an approximate circular shape. We attempted to make this process realistic to real-world puzzle solving by not being too precise about getting the pieces as close together as possible. In effect, arranging the pieces too carefully in a grid would impact the assumption of puzzle piece circularity. We measured the length $X_s$ and width $Y_s$ of the resulting oval of pieces (Figure \ref{fig:puzzle_measurement}) and calculated the area using the formula for the area of an ellipse: $A = \frac{\pi}{4} X_s Y_s$. For the two 1000-piece puzzles (last two rows of Table \ref{tab:data}), we laid out the pieces in an approximate rectangular shape instead of an oval and calculated the area of a rectangle instead of an ellipse.

\begin{figure}
\centering
\includegraphics[width=0.7\linewidth]{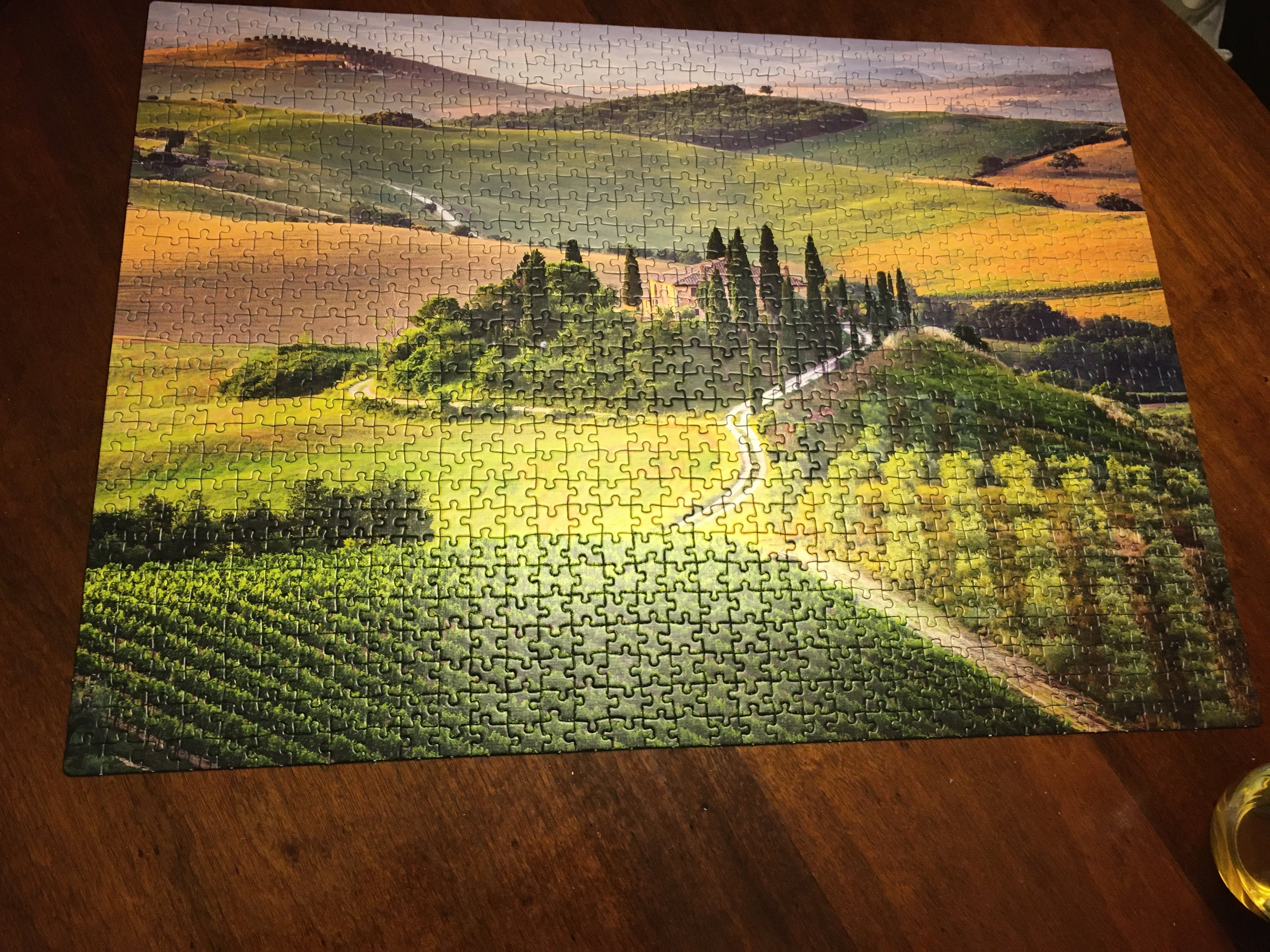} 
\caption{Photograph of an assembled 1008-piece puzzle. We measured the length and width ($X_a$ and $Y_a$) of the assembled puzzle to estimate the assembled area.} \label{fig:puzzle_asssembled}
\end{figure}

\begin{table} 
\begin{center}
\begin{tabular}{ |c|c|c|c|c|} 
 \hline
 $N$ & $X_a$ (cm)  & 	$Y_a$ (cm) &	$X_s$ (cm) &	$Y_s$ (cm) \\
 \hline 
	1008 &	50.2 &		69.0 &	83.0 &	85.0 \\
	252 &	26.6 &		34.4 &	45.0 &	46.5  \\
 	9 &		15.6 &		21.4 &	25.9 &	23.3  \\
 	500 &	50.8 &		50.9 &	78.4 &	74.8  \\
 	1026 &	67.8 &		48.9 &	88.0 &	86.8  \\
 	27 &	296.8 &		14.6 &	107.2&	92.0  \\
 	2000 &	99.1 &		68.6 &	123.1 &	131.0 \\
 	1000 &	50.8 &		68.5 &	112.0 &	69.0  \\
    1000 &	99.3 &		33.0 &	132.4 &	57.5  \\
 \hline
\end{tabular} 
\end{center}
\caption{Measurements of puzzle unassembled and assembled dimensions for 9 puzzles. We estimated the uncertainty in our measurement of $X_a$ and $Y_a$ as $0.2$ cm and the uncertainty in in $X_s$ and $Y_s$ as $0.5$ cm.}
\label{tab:data}
\end{table}

\begin{figure}
\centering
\includegraphics[width=1\linewidth]{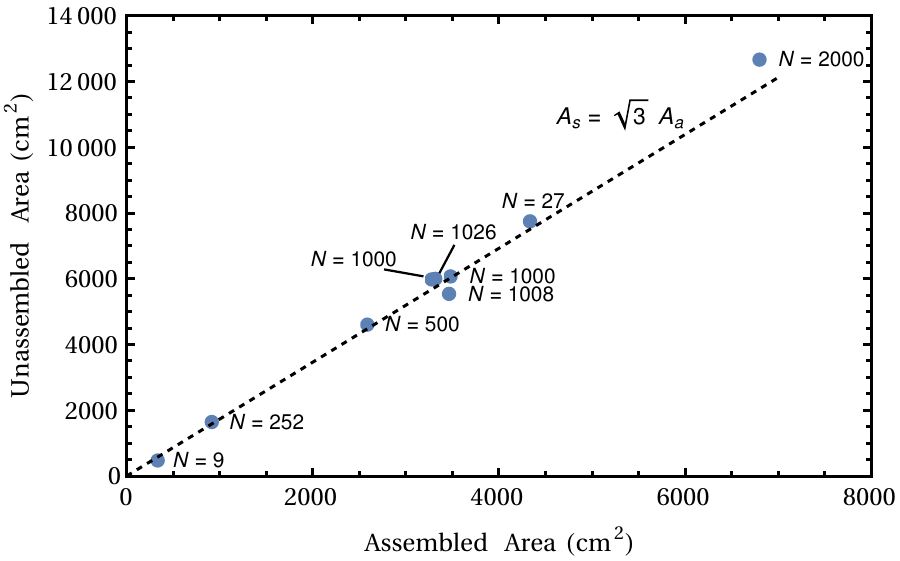} 
\caption{\textbf{Puzzle unassembled area is accurately predicted by theory that depends only on the assembled area.} Measured puzzle unassembled areas are plotted as a function of assembled areas (blue dots). Error bars are obscured by dot size. The theoretical predition $A_s = \sqrt{3}A_a$ is plotted as a black dashed line. No fitting was performed.} \label{fig:results}
\end{figure}

We assembled each puzzle (Figure \ref{fig:puzzle_asssembled}), then measured the lengths and widths of the assembled puzzles (Table \ref{tab:data}). We calculated the corresponding areas and plotted the results alongside the theoretical prediction (Figure \ref{fig:results}). We found close agreement between realistic measurements and our theoretical prediction across a wide range of puzzle areas and numbers of pieces.  

\section{Discussion}

In this work, we derive a simple formula based on physical principles for the area of an unassembled puzzle ($A_s$) based on its assembled area ($A_a$): $A_s = \sqrt{3} A_a$. For a puzzle of area $A_a$, the assembler will need a surface slightly less than twice that area to lay out all the pieces flat.

Surprisingly, this result does not depend on the number of pieces. This is counterintuitive, but can be understood with physical reasoning based on the constitutive equations. Consider two puzzles with identical assembled areas but different numbers of pieces. With a small number of large pieces, the gaps between pieces are larger (the hexagonal lattice spacing is larger), but this area is multiplied by a small number of pieces. Conversely, for many small pieces, the lattice spacing is smaller but there are more pieces. Each lattice element is multiplied by $N$, but there is a factor of $1/N$ in the area of each lattice element; these cancel out, giving a final result independent of $N$.

\bibliographystyle{unsrt}
\bibliography{puzzle_refs.bib}

\begin{thebibliography}{1}

\bibitem{Ashcroft76}
N.~W. Ashcroft and N.~D. Mermin.
\newblock {\em {S}olid {S}tate {P}hysics}.
\newblock Holt-Saunders, 1976.

\bibitem{Chang2010}
Hai-Chau Chang and Lih-Chung Wang.
\newblock {A Simple Proof of Thue's Theorem on Circle Packing}.
\newblock {\em arXiv:1009.4322}, 2010.

\end{thebibliography}

\end{document}